\documentclass[12pt]{article}
\usepackage{amssymb,amsmath,amsfonts,bbm,bbold,accents,enumitem} 
\author{Fernando Argentieri }

\font\teneufm=eufm10
\font\seveneufm=eufm7
\font\fiveeufm=eufm5
\newfam\eufmfam
\textfont\eufmfam=\teneufm
\scriptfont\eufmfam=\seveneufm
\scriptscriptfont\eufmfam=\fiveeufm

\newcommand\beq[1]{ \begin{equation}\label{#1} }
\newcommand{\eeq}{ \end{equation} }

\newcommand\beqa[1]{ \begin{eqnarray} \label{#1}}
\newcommand{\eeqa}{ \end{eqnarray} }
\newcommand{\beqano}{ \begin{eqnarray*} }
\newcommand{\eeqano}{ \end{eqnarray*} }

\newtheorem{theorem}{Theorem}
\newtheorem{definition}{Definition}
\newtheorem{proposition}{Proposition}
\newtheorem{lemma}{Lemma}
\newtheorem{sublemma}{Sublemma}
\newtheorem{remark}{Remark}
\newtheorem{notationalremark}{Notations}
\newtheorem{corollary}{Corollary}
\newtheorem{assumption}{Assumption}
\newtheorem{claim}{Claim}
\newtheorem{tools}{$\negsp\negsp$}[subsection]

\newcommand{\II}{\mathcal{I}}
\newcommand\thm[1]{ \begin{theorem}\label{#1}}
\newcommand\thmtwo[2]{ \begin{theorem}[#1]\label{#2}}
\newcommand\ethm{ \end{theorem} }
\newcommand\dfn[1]{ \begin{definition}\label{#1} \rm}
\newcommand\dfntwo[2]{ \begin{definition}[#1]\label{#2} \rm}
\newcommand\edfn{ \end{definition} }
\newcommand\pro[1]{ \begin{proposition}\label{#1}}
\newcommand\protwo[2]{ \begin{proposition}[#1]\label{#2}}
\newcommand\epro{ \end{proposition} }
\newcommand\lem[1]{ \begin{lemma}\label{#1}}
\newcommand\lemtwo[2]{ \begin{lemma}[#1]\label{#2}}
\newcommand\elem{ \end{lemma} }
\newcommand\sublem[1]{ \begin{sublemma}\label{#1}}
\newcommand\sublemtwo[2]{ \begin{sublemma}[#1]\label{#2}}
\newcommand\esublem{ \end{sublemma} }
\newcommand\rem[1]{ \begin{remark}\label{#1} \rm}
\newcommand\erem{ \end{remark} }
\newcommand\notrem[1]{ \begin{notationalremark}\label{#1} \rm}
\newcommand\enotrem{ \end{notationalremark} }
\newcommand\cor[1]{ \begin{corollary}\label{#1}}
\newcommand\cortwo[2]{ \begin{corollary}[#1]\label{#2}}
\newcommand\ecor{ \end{corollary} }
\newcommand\asmp[1]{ \begin{assumption}\label{#1}}
\newcommand\asmptwo[2]{ \begin{assumption}[#1]\label{#2}}
\newcommand\easmp{ \end{assumption} }
\newcommand\clm[1]{ \begin{claim}\label{#1}}
\newcommand\eclm{ \end{claim} }
\newcommand{\proof}{\par\medskip\noindent{\bf Proof\ }}
\expandafter\chardef\csname pre amssym.def
at\endcsname=\the\catcode`\@
\catcode`\@=11
\def\undefine#1{\let#1\undefined}
\def\newsymbol#1#2#3#4#5{\let\next@\relax
 \ifnum#2=\@ne\let\next@\msafam@\else
 \ifnum#2=\tw@\let\next@\msbfam@\fi\fi
 \mathchardef#1="#3\next@#4#5}
\def\mathhexbox@#1#2#3{\relax
 \ifmmode\mathpalette{}{\m@th\mathchar"#1#2#3}%
 \else\leavevmode\hbox{$\m@th\mathchar"#1#2#3$}\fi}
\def\hexnumber@#1{\ifcase#1 0\or 1\or 2\or 3\or 4\or 5\or 6\or 7\or
8\or
 9\or A\or B\or C\or D\or E\or F\fi}
\ifcase\@ptsize
 \font\tenmsb=msbm10
 \font\sevenmsb=msbm7
 \font\fivemsb=msbm5
\or
 \font\tenmsb=msbm10 scaled \magstephalf
 \font\sevenmsb=msbm7 scaled \magstephalf
 \font\fivemsb=msbm5  scaled \magstephalf
\or
 \font\tenmsb=msbm10 scaled \magstep1
 \font\sevenmsb=msbm7 scaled \magstep1
 \font\fivemsb=msbm5 scaled \magstep1
\fi
\newfam\msbfam
\textfont\msbfam=\tenmsb
\scriptfont\msbfam=\sevenmsb
\scriptscriptfont\msbfam=\fivemsb
\edef\msbfam@{\hexnumber@\msbfam}
\def\Bbb#1{\fam\msbfam\relax#1}
\def\widehat#1{\setboxz@h{$\m@th#1$}%
 \ifdim\wdz@>\tw@ em\mathaccent"0\msbfam@5B{#1}%
 \else\mathaccent"0362{#1}\fi}
\def\widetilde#1{\setboxz@h{$\m@th#1$}%
 \ifdim\wdz@>\tw@ em\mathaccent"0\msbfam@5D{#1}%
 \else\mathaccent"0365{#1}\fi}

\def\RIfM@{\relax\ifmmode}
\def\nonmatherr@#1{\errmessage{\string#1\space allowed only in math mode}}
\def\Bbb{\RIfM@\expandafter\Bbb@\else
 \expandafter\nonmatherr@\expandafter\Bbb\fi}
\def\Bbb@#1{{\Bbb@@{#1}}}
\def\Bbb@@#1{\fam\msbfam\relax#1}
\def\setboxz@h{\setbox\z@\hbox}
\def\wdz@{\wd\z@}
\catcode`\@=\csname pre amssym.def at\endcsname
\newcommand{\acapo}{\linebreak}

\newcommand{\Giu}{{\bigskip\noindent}}
\newcommand{\nl}{{\smallskip\noindent}}

\newcommand{\qed}{\hskip.5truecm
\vrule width 1.7truemm height 3.5truemm depth 0.truemm
\par\Giu}

\newcommand{\negsp}{\hspace{-.09truecm}}  

\newcommand{\ti}{{\mathcal{I}}^{1}_{\t}}
\newcommand{\tiii}{{\mathcal{I}}^{3}_{\t}}
\newcommand{\tii}{{\mathcal{I}}^{2}_{\t}}
\newcommand{\igt}{{\mathcal{I}}^{1}_{\g,\t}}
\newcommand{\iigt}{{\mathcal{I}}^{2}_{\g,\t}}
\newcommand{\iiigt}{{\mathcal{I}}^{3}_{\g,\t}}
\newcommand{\IIt}{\mathcal{I}_{\t}}
\newcommand{\DDt}{\mathcal{D}_{\t}}
\newcommand{\qmg}{\frac{\g}{q_{m}^{\t+1}}}
\newcommand{\qm}{\frac{p_{m}}{q_{m}}}
\newcommand{\pq}{\frac{p}{q}}
\newcommand{\qnnng}{\frac{\g}{q_{n+2}^{\t+1}}}
\newcommand{\qng}{\frac{\g}{q_{n}^{\t+1}}}

\newcommand{\xn}{\frac{1}{\g}-\frac{q_{n+1}}{q_{n}^{\t}}}
\newcommand{\nppp}{\frac{p_{n+2}}{q_{n+2}}}
\newcommand{\npp}{\frac{p_{n+1}}{q_{n+1}}}
\newcommand{\np}{\frac{p_{n}}{q_{n}}}
\newcommand{\Dgt}{D_{\gamma,\tau}}
\newcommand{\gta}{\gamma(\alpha,\tau)}
\newcommand{\gnat}{\g_{n}(\a,\t)}
\newcommand{\gtal}{\g_{-}(\a,\t)}

\newcommand{\gtar}{\g_{+}(\a,\t)}

\newcommand{\gmat}{\g_{m}(\a,\t)}

\newcommand{\dst}{\displaystyle}

\newcommand\su[1]{ \frac{1}{ {#1}} }

\renewcommand{\a }{ {\alpha}   }
\renewcommand{\b}{ {\beta}   }
\newcommand{\g}{ {\gamma}   }

\newcommand{\e }{ {\epsilon}   }

\newcommand{\m}{ {\mu}   }

\renewcommand{\t}{ {\tau}   }

\newcommand{\implica}{\ \Longrightarrow\ }

\renewcommand\subset{\subseteq}

\title{Diophantine sets in general are Cantor sets}
\author{Fernando Argentieri}
\begin{document}
\maketitle

\begin{abstract}
 Let $\g\in(0;\su{2}),\t\geq 1$ and define the ``$\g,\t$ Diophantine set" as: 
$$\Dgt:=\{\a\in (0;1): ||q\a||\geq\frac{\g}{q^{\t}}\quad\forall q\in\Bbb{N}\},\qquad||x||:=\inf_{p\in\Bbb{Z}}|x-p|. $$ 
In this paper we study the topology of these sets and we show that, for large $\t$ and for almost all $\g>0$, $\Dgt$ is a Cantor set.

\end{abstract}
\section{Introduction}
Diophantine sets play an important role in dynamical systems, in particular, in
small divisors problems with applications to KAM theory, Aubry-Mather theory,
conjugation of circle diffeomorphisms, etc. (see, for example, [3], [5], [9], [12], [13], [14], [16]). 

\nl
The set $\Dgt$ is compact and totally
disconnected (since $\Dgt\cap\Bbb{Q}=\emptyset$), however, these sets may be not Cantor sets. In fact, in \cite{17} we have shown various examples in which $\Dgt$ have isolated points. In this paper we prove the following:\acapo
{\bf Theorem}
Let $\t>\frac{3+\sqrt{17}}{2}$. Then, for almost all $\g>0$ $\Dgt$ is a Cantor set.

\nl
By \cite{6}, for $\t=1$ and $\su{3}<\g<\su{2}$ $\Dgt$ is countable (and non empty for $\g>\su{3}$ small enough). In particular, this result does not holds for $\t=1$. We expect that $\t>\frac{3+\sqrt{17}}{2}$ can be improved with $\t>3$. However it is not clear what is the best constant. Following the same proof, we can prove also that, fixed $\t>\frac{3+\sqrt{17}}{2}$, for almost all $\g>0$, if $\a\in\Dgt$ and $U$ is an open neighborhood that contains $\a$, then $\m(\Dgt\cap U)>0$.\acapo

\nl
 The paper is organized as follows: in the second section we give some basic definitions and remarks, in the third section we prove our result and, in the last section are present some natural questions.

\section{Definitions and remarks}

\subsection{Definitions}

\begin{itemize}
    \item $\Bbb{N}:=\{1,2,3,...\}$, $\Bbb{N}_{0}:=\{0,1,2,3,...\}$ 
    \item Given $a,b\in\Bbb{Z} -\{0\}$, we indicate with $(a,b)$ the maximum common divisor of $a$ and $b$.
    \item Let $\a$ be a real number. We indicate with $[\a]$ the integral part of $\a$, with $\{\a\}$ the fractional part of $\a$ .
     \item Given E$\subset{\Bbb{R}}$, we indicate with $\mathcal{I}$(E) the set of isolated points of E.
     \item Given E$\subset{\Bbb{R}}$, we indicate with $\mathcal{A}$(E) the set of accumulated points of E.
     \item We say that E$\subset{\Bbb{R}}$ is perfect if $\mathcal{A}$(E)=E.
    \item Given a Borel set E$\subset{\Bbb{R}}$ we denote with $\m$(E) the Lebesgue measure of E.
    \item A topological space X is a totally disconnected space if the points are the only connected subsets of X.
    \item $X\subset\Bbb{R}$ is a Cantor set if it is closed, totally disconnected and perfect.
    \item For $E\subset\Bbb{R}^{n}$, $\dim_{H}E$ is the Hausdorff dimension of $E$.
    \item Given $\a\in\Bbb{R}$ we define:
    $$||\a||:=\min_{p\in\Bbb{Z}}|\a -p|$$
    \item Given $\g>0, \t\geq1$, we define the $(\g,\t)$ Diophantine points in $(0;1)$ as the numbers in the set:
    $$ \Dgt:=\{\a\in (0 ;1): ||q\a||\geq\frac{\g}{q^{\t}}\quad\forall q\in\Bbb{N}\}$$
    \item $$ D^{\Bbb{R}}_{\g,\t}:=\{\a\in\Bbb{R}:||q\a||\geq\frac{\g}{q^{\tau}}\quad\forall q\in\Bbb{N}\},$$ $$D_{\t}:=\bigcup_{\g>0}D_{\g,\t},\quad D:=\bigcup_{\t\geq 1}D_{\t}.$$  We call $D$ the set of Diophantine numbers.
    \item Given $\t\geq 1,\a\in \Bbb{R}$, we define:
    $$\gta:=\inf_{q\in\Bbb{N}}q^{\t}||q\a||$$
    \item Given $\a\in\Bbb{R}$ we define:
    $$\t(\a):=\inf\{\t\geq 1:\gta>0\}$$
   \item Given an irrational number $\a=[a_{0};a_{1},...]:=a_0+\su{a_1+\su{a_2+...}}$, we denote with $\{\frac{p_n}{q_n}\}_{n\in\Bbb{N}_{0}}$ the convergents of $\a$, $\a_{n}:=[a_{n};a_{n+1},...]$\footnote{for informations about continued fractions see [4],[8],[15] }. 
   \item We indicate with $[a_1,a_2,a_3,...]:=\su{a_1+\su{a_2+\su{a_3+...}}}$.
   \item Let $\a$ be an irrational number. We define:
   $$\gnat:=q_n^{\t}||q_n\a||=q_n^{\t}|q_n\a-p_n|$$
   \item Let $\t\geq 1$,
   $$\gtal:=\inf_{n\in2 \Bbb{N}_{0}}\gnat,$$
   $$\gtar:=\inf_{n\in2\Bbb{N}_{0}+1}\gnat,$$
   $${\DDt}:=\{\a\in D_{\t} :\t(\a)=\t\},$$
   $$\igt:=\{\a\in D_{\g,\t}: \exists n\not\equiv m\quad{(\rm{mod} 2)}, \gnat=\gmat=\gta\},$$
   $$\iigt:=\{\a\in D_{\g,\t}: \exists n\in\Bbb{N}_{0}, \gnat=\g(\a,\t)\}\cap(\igt)^{c},$$
   $$\iiigt:={\mathcal{I}}(\Dgt)\cap(\igt\cup\iigt)^{c},$$
   $$\ti:=\bigcup_{\g>0}\igt,$$
   $$\tii:=\bigcup_{\g>0}\iigt,$$
   $$\tiii:=\bigcup_{\g>0}\iiigt.$$
   \end{itemize}
\subsection{Remarks}
We list here some simple remarks. For a proof see \cite{17}.
\begin{enumerate}
    \item $\a\in \Dgt\iff 1-\a\in\Dgt$.
    \item $\gta\leq \min\{\a,1-\a\}.$
    \item Fixed $\t\geq 1$, $\g(.,\t):D_{\t} \rightarrow (0,\frac{1}{2})$.
    \item $\Dgt^{\Bbb{R}}=\bigcup_{n\in\Bbb{Z}}(\Dgt+n)$, thus we can restrict to study the Diophantine points in $(0,1)$.
    \item 
    \beq{fond}
    \left\{
    \begin{array}{l}
           \gnat=\frac{q_{n}^{\t}}{\a_{n+1}q_{n}+q_{n-1}},\\
           \su{\gnat}=\frac{q_{n+1}}{q_{n}^{\t}}+\frac{1}{\a_{n+2}q_{n}^{\t-1}}
    \end{array}\right.
    \eeq
    \item $\gta=\inf_{n\in\Bbb{N}_{0}}\gnat$.
    \item If $\t<\t(\a)$, then $\gta=0$; if $\t>\t(\a)$ then $\gta>0$. Moreover, for $\t>\t(\a)$ the inf is a minimum.
    \item $\a\in {\DDt}\iff \t(\a)=\t$ and $\gta>0$.
    \item If $\a\in\igt$, then $\a$ is an isolated point of $\Dgt$.
    \item The cardinality of ${\ti}$ is at most countable.
    \item $\m({\DDt})=0$ for all $\t\geq 1$.
    \item $\g_{0}(\a,\t)=\left\{\a\right\}$, in particular $\g_{0}(\a,\t)$ does not depend on $\t$.
    \item Let $\pq$ a rational number. 
    \beq{}
    \a\in D_{\t}\iff \left\{\a+\pq\right\}\in D_\t,
    \eeq
    \beq{}
    \a\in{\DDt}\iff \left\{\a+\pq\right\}\in{\DDt}.
    \eeq
    \item If $\t>\t(\a)$, $\gtal=\gtar$, then $\a\in{\IIt}$.
    \item $\a\in D_\t\iff q_{n+1}=O(q_{n}^{\t}).$
\end{enumerate}

\section{Proof of the Theorem}
In the first part of this section, we suppose without loss of generality, that $n$ is always even. In fact, for $n$ odd it suffices to consider $1-\a$ ($\a=[a_{0};...,a_{n},...]\in \Dgt\iff 1-\a\in\Dgt$, and the denominators of the odd convergents to $1-\a$ are the same of the even convergents to $\a$, hence, by symmetry, all that is demonstrated for $n$ even continues to hold if $n$ is odd). Moreover, in all the section $0<\g<\su{2}$ (otherwise $\Dgt=\emptyset$).
We want to prove that, for $\t>\frac{3+\sqrt{17}}{2}$:
$$\m\left(\left\{0<\g<\frac{1}{2}:{\II}(\Dgt)\not=\emptyset\right\}\right)=0.$$
By Remark (j) it is enough to prove it for $\iigt$ and $\iiigt$.
Observe that the isolated points of type 2,3 are obtained by infinitely many intersections of intervals centered in rational numbers $\pq$ with length $\frac{2\g}{q^{\t+1}}$. Thus, the first step is to show that, given $\a\in\Dgt$, it is enough (up to a set of measure zero and for $\t$ big enough) to control the intersection of intervals centred in the convergents. The second step will be to show that, if intervals centred in the convergents intersects, then the coefficients of the continued fractions cannot grow too. In the final step we prove that, when intervals centred in the convergents do not intersect and for big convergents, the interval between two subsequent convergentes (with the same parity) contains a diophantine sets with positive mesure.

\lem{}
\label{l8}
Let $\g>0, \t>1,\a\in\Dgt$, $\np$ the convergents to $\a$, $$I_{n}:=\left(\np,\nppp\right).$$ Suppose that $\exists N\in\Bbb{N}$ such that, for all $n>N$ even:
\beq{2}
\np+\qng<\nppp-\qnnng.
\eeq
For $n>N$ define $$A_{n}:=\left(\np+\qng, \nppp-\qnnng\right).$$
Moreover, suppose that for every $n$ (even):
\beq{1}
\a-\np>\qng
\eeq
Then, there exists $ N_{1}\in\Bbb{N}$ such that, for all $n>N_{1}$:
$$\frac{p}{q}\not\in I_{n}\implica\frac{p}{q}+\frac{\g}{q^{\t+1}},\frac{p}{q}-\frac{\g}{q^{\t+1}}\not\in A_{n}.$$
\elem{}
\proof
Note that it is enough to verify the inequality when $\pq<\a$. In fact the inequality is trivial if $\pq>\a$ (because of $\a\in\Dgt$ implies $\pq-\frac{\g}{q^{\t+1}}\geq\a>\frac{p_{n+2}}{q_{n+2}}+\frac{\g}{q_{n+2}^{\t+1}}$ by (\ref{1})). 
By (\ref{2}) it follows that $A_{n}\cap A_{m}=\emptyset$ for $n\not=m$, with $n,m>N$ even. From $$\a-\np>\qng$$ for $n$ even, we get $$\max_{2n\leq N} \frac{p_{2n}}{q_{2n}}+\frac{\g}{q_{2n}^{\t+1}}=:C<\a,$$ from which it follows that there exists $N_{1}\in\Bbb{N}$ such that for $n$ even, $n>N_{1}$: $$\np-\qng>C.$$
If $\frac{p}{q}=\frac{p_{m}}{q_{m}}\not \in I_{n}$ is an even convergent to $\a$ with $n>N_{2}:=\max\{N,N_{1}\}$ then, for $m\leq N$ even: $$\frac{p_{m}}{q_{m}}<\np.$$ 
Moreover, by definition of $N_1$ it follows that:
$$\frac{p_{m}}{q_{m}}+\frac{\g}{q_{m}^{\t+1}}\leq C<\np-\qng,$$
from which it follows that the Lemma holds if $\frac{p}{q}=\frac{p_{m}}{q_{m}}$ is an even convergent to $\a$ with $m\leq N$. If $m>N$ and $n>m$ is even:
$$\qm+\qmg<\frac{p_{m+2}}{q_{m+2}}-\frac{\g}{q_{m+2}^{\t+1}}\leq \np+\qng$$
while, for $n<m$ even:
$$\qm-\qmg>\frac{p_{m-2}}{q_{m-2}}+\frac{\g}{q_{m-2}^{\t+1}}\geq \nppp+\qnnng.$$
So Lemma \ref{l8} is true if $\pq$ is an even convergent to $\a$.
Thus, Lemma \ref{l8} remains to be verified  when $\pq$ is not a convergent to $\a$. It is no restrictive to suppose that there exists $m\not=n$ even for which $\frac{p}{q}\in I_{m}$, otherwise Lemma \ref{l8} is trivial. Now we show that, for $m$ big enough:
$$\frac{p}{q}+\frac{\g}{q^{\t+1}}, \frac{p}{q}-\frac{\g}{q^{\t+1}}\in\left(\frac{p_{m}}{q_{m}}-\frac{\g}{q_{m}^{\t+1}},\frac{p_{m+2}}{q_{m+2}}+\frac{\g}{q_{m+2}^{\t+1}}\right)$$
from which Lemma \ref{l8} follows immediately by (\ref{1}). By the properties of Farey sequence, for the rationals $\frac{p}{q}\in I_{m}$ we have $q>q_{m}$, so the inequality: $$\frac{p}{q}-\frac{\g}{q^{\t+1}}>\frac{p_{m}}{q_{m}}-\frac{\g}{q_{m}^{\t+1}}$$ holds. It remains to show that:
$$\frac{p}{q}+\frac{\g}{q^{\t+1}}<\frac{p_{m+2}}{q_{m+2}}+\frac{\g}{q_{m+2}^{\t+1}}.$$
This inequality holds for $q\geq \frac{q_{m+2}}{2}$ and $m$ big enough. In fact, in that case:
$$\frac{p_{m+2}}{q_{m+2}}-\pq\geq \frac{1}{q q_{m+2}}> \frac{\g}{q^{\t+1}}-\frac{\g}{q_{m+2}^{\t+1}},$$
that is true for $m$ big enough (because of $\t>1$). So, we can assume that $q_{m}<q<\frac{q_{m+2}}{2}$. Because we have assumed that  $\pq$ is not a convergent, by Legendre's Theorem (see \cite{8}), we have: $$\a-\pq>\frac{1}{2 q^{2}},$$
while, because $\frac{p_{m}}{q_{m}}$ is a convergent, we have:
$$\a-\frac{p_{m+2}}{q_{m+2}}<\frac{1}{q_{m+2}^{2}}.$$
So, putting together the two inequalities, if $q<\frac{q_{m+2}}{2}$:
$$\frac{p_{m+2}}{q_{m+2}}-\frac{p}{q}=\frac{p_{m+2}}{q_{m+2}}-\a+\a-\frac{p}{q}> \frac{1}{2q^{2}}-\frac{1}{q_{m+2}^{2}}>-\frac{\g}{q_{m+2}^{\t+1}}+\frac{\g}{q^{\t+1}}\iff$$
$$\frac{1}{2q^{2}}-\frac{\g}{q^{\t+1}}>\frac{1}{q_{m+2}^{2}}-\frac{\g}{q_{m+2}^{\t+1}},$$
that is true for $m$ big enough (it follows by $q_{m}<q<\frac{q_{m+2}}{2}$). So Lemma \ref{l8} is proved. \qed
\\
\nl
We know by Farey's sequence that for $\pq\in I_n$, $q>q_{n+1}$. So, there are a finite numbers of $\pq\in I_n$ with $q<q_{n+2}$. In the next Lemma we want to control the distance between these numbers and $\nppp-\frac{\g}{q_{n+2}^{\t+1}}$.

\lem{}
\label{l9}
Let $\g>0$, $\t>3, \a\in\Dgt, \np$ the convergents to $\a$. 
There exists $N_{1}\in\Bbb{N}$ such that, for $n>N_{1}$:
$$\frac{p}{q}\in I_{n}, q<q_{n+2}\implica \frac{p}{q}+\frac{\g}{q^{\t+1}}<\nppp-\qnnng -\frac{2\g}{q_{n+2}^{\t-1}}.$$
\elem{}
\proof
Let $n>N$, $\frac{p}{q}\in I_{n}$, so by definition of convergents and the fact that $\np<\frac{p}{q}<\nppp$ we get that $\frac{p}{q}$ is not a convergent.
If $q\geq\frac{q_{n+2}}{2}$ we get:
$$\nppp-\pq\geq\frac{1}{q q_{n+2}}\geq \frac{1}{q_{n+2}^{2}}>\frac{\g 2^{\t+1}}{q_{n+2}^{\t+1}}+\qnnng+\frac{2\g}{q_{n+2}^{\t-1}}\geq \frac{\g}{q^{\t+1}}+\qnnng+\frac{2\g}{q_{n+2}^{\t-1}}$$
for $n$ big enough (because of $\t>3$).
So, for $n$ big enough, the inequality remain to be proved for $q<\frac{q_{n+2}}{2}$.
In that case:
$$\nppp-\frac{p}{q}=\nppp-\a+\a-\pq>\frac{1}{2q^{2}}-\frac{1}{q_{n+2}^{2}}>\frac{\g}{q^{\t+1}}+\qnnng+\frac{2\g}{q_{n+2}^{\t-1}}\iff$$
$$\frac{1}{2q^{2}}-\frac{\g}{q^{\t+1}}>\frac{1}{q_{n+2}^{2}}+\qnnng+\frac{2\g}{q_{n+2}^{\t-1}}.$$
From the fact that $$G(x):=\frac{1}{2x^{2}}-\frac{\g}{x^{\t+1}}$$ is a decreasing function for $x$ big enough, it is enough to show the inequality for $q=[\frac{q_{n+2}}{2}]$. In this case we get:

$$\frac{1}{2q^{2}}-\frac{1}{q_{n+2}^{2}}\geq\frac{2}{q_{n+2}^{2}}-\frac{1}{q_{n+2}^{2}}=\frac{1}{q_{n+2}^{2}}>\frac{\g}{q^{\t+1}}+\qnnng+\frac{2\g}{q_{n+2}^{\t-1}}$$ for $n$ big enough (for $\t>3$), so $\exists N_{1}\in\Bbb{N}$ such that, when $n>N_{1}$ is even the inequality is verified. \qed
\\
\nl
\lem{}
\label{l10}
Let $\t>3$ $,\a=[a_{1},a_{2},...]\in\Dgt, \np$ the convergents to $\a$, then $\exists N\in\Bbb{N}$ such that for all $n>N$ even:
$$\m\left(\bigcup_{\frac{p}{q}\in I_{n}, q\geq q_{n+2}}\left(\frac{p}{q}-\frac{\g}{q^{\t+1}}, \frac{p}{q}+\frac{\g}{q^{\t+1}}\right)\right)<\frac{2\g}{q_{n+2}^{\t-1}}$$
\elem{}
\proof
$$ \m\left(\bigcup_{\frac{p}{q}\in I_{n}, q\geq q_{n+2}}\left(\frac{p}{q}-\frac{\g}{q^{\t+1}}, \frac{p}{q}+\frac{\g}{q^{\t+1}}\right)\right)$$
$$<\sum_{q\geq q_{n+2}}\sum_{q\np<p<q\nppp}\frac{2\g}{q^{\t+1}}<2\g\left(\nppp-\np\right)\sum_{q\geq q_{n+2}}\su{q^{\t}}$$
$$<2\g C\left(\nppp-\np\right)\su{q_{n+2}^{\t-1}}=o\left(\frac{2\g }{q_{n+2}^{\t-1}}\right)$$
for some constant $C>0$.\qed
\\
\nl
\lem{}
\label{l11}
Let $\t>1, \g>0,$ $\a=[a_{1},a_{2},...]\in\Dgt, \np$ be the convergents to $\a$. Then:
\beq{*}
\quad \np+\qng<\nppp-\qnnng\iff
\eeq
\beq{}
a_{n+2}>\frac{q_{n}}{\g q_{n+1}}\su{(\xn)-\frac{q_{n}q_{n+1}}{ q_{n+2}^{\t+1}}}-\frac{q_{n}}{q_{n+1}} \quad
\eeq
\elem{}

\proof
(\ref{*}) is true if and only if:
$$\nppp-\np=\nppp-\npp+\npp-\np=$$
$$\frac{1}{q_{n}q_{n+1}}-\frac{1}{q_{n+1}q_{n+2}}>\qnnng+\qng\iff$$
$$\frac{1}{q_{n+2}q_{n+1}}<\frac{1}{q_{n}q_{n+1}}-\qng-\qnnng\iff$$
$$\frac{1}{q_{n+2}q_{n+1}}<\frac{\g}{q_{n}q_{n+1}}(\xn)-\qnnng\iff$$
$$\su{q_{n+2}}<\frac{\g}{q_{n}}(\xn)-q_{n+1}\qnnng\iff $$
\beq{}
\left\{
\begin{array}{l}
\dst\xn>\frac{q_{n}q_{n+1}}{ q_{n+2}^{\t+1}},\\ \ \\
\dst q_{n+2}>\frac{q_{n}}{\g}\su{(\xn)-\frac{q_{n}q_{n+1}}{ q_{n+2}^{\t+1}}} \\ 
\end{array}\right.
\ \\\eeq{}
The first inequality is always true because of:
$$\xn>\frac{1}{\a_{n+2}q_{n}^{\t-1}}>\frac{q_{n}q_{n+1}}{ q_{n+2}^{\t+1}}.$$

\Giu
So Lemma \ref{l11} follows from the fact that $q_{n+2}=a_{n+2}q_{n+1}+q_{n}$.\qed

\lem{}
\label{l12}
Let $\t>1$, for almost all $\g\in(0,\su{2})$ (for $\g\geq\su{2}$ $\Dgt=\emptyset$), given $\e>0$ there exists $C=C(\e,\g)>0$ such that:
$$\left|\su{\g}-\frac{p}{q^{\t}}\right|\geq \frac{C}{q^{\t+1+\e}}$$
for all $\frac{p}{q}\in\Bbb{Q}$.
\elem{}
\proof
Define $B_{C,k}:=\left\{\a:|\a-\frac{p}{q^{\t}}|\geq \frac{C}{q^{k}}\quad\forall \frac{p}{q}\in\Bbb{Q}\right\}$,
so $\a\in B_{C,k}^{c}\iff$ there exists $\frac{p}{q}$ such that $\a\in \left(\pq-\frac{C}{q^{k}},\pq+\frac{C}{q^{k}}\right)$. So, given $N\in\Bbb{N}$ we get:
$$\m\left(B_{C,k}^{c}\cap \left(-N,N\right)\right)<\sum_{q>0}\sum_{-N q^{\t}<p<N q^{\t}} \frac{2C}{q^{k}}<\sum_{q>0}\frac{4NC}{q^{k-\t}}$$
and for $k>\t+1$, $C$ that tends to zero, also $$\m\left(B_{C,k}^{c}\cap \left(-N,N\right)\right)$$ goes to zero. From the arbitrariness of $N$ we obtain: $$\m\left(\bigcap_{C>0} B_{C,k}^{c}\right)=0$$ for $k>\t+1$, from which follows Lemma \ref{l12}. \qed

\nl
\lem{}
\label{l14}
Let $\t>1$, $\a=[a_{1},a_{2},...]\in\Dgt$, $\np$ the convergents to $\a$. The inequality:
\beq{**}
\quad\np+\qng<\nppp-\qnnng-\frac{2\g}{q_{n+2}^{\t-1}}
\eeq
 is definitively verified if and only if definitively:
\beq{ci}
{a_{n+2}>\frac{q_{n}}{\g q_{n+1}}\su{(\xn)-\frac{q_{n}q_{n+1}}{q_{n+2}^{\t+1}}-\frac{2q_{n}q_{n+1}}{q_{n+2}^{\t-1}}}-\frac{q_{n}}{q_{n+1}}}
\eeq

\nl
\elem{}
\rem{}
\label{r16}
Observe that (\ref{ci}) is definitively true if: $$\limsup\frac{q_{n+1}}{q_{n}^{\t}}<\frac{1}{\g},$$ because in that case:
$$\limsup \frac{q_{n}}{\g q_{n+1}}\su{(\xn)-\frac{q_{n}q_{n+1}}{ q_{n+2}^{\t+1}}-\frac{2q_{n}q_{n+1}}{q_{n+2}^{\t-1}}}-\frac{q_{n}}{q_{n+1}}<1.$$
Thus, if for infinitely many $n$ even $(\ref{ci})$ is not verified, for this $n$, with $n$ big enough: $$\frac{q_{n+1}}{q_{n}^{\t}}\sim \frac{1}{\g},$$ so $q_{n+1}\sim \frac{q_{n}^{\t}}{\g}$.
\erem{}

\proof 

In a similar way of Lemma \ref{l11}, (\ref{**}) is verified if and only if:
\beq{}
\left\{
\begin{array}{l}
\dst \xn>\frac{q_{n}q_{n+1}}{ q_{n+2}^{\t+1}}+\frac{2q_{n}q_{n+1}}{q_{n+2}^{\t-1}},\\  \ \\
\dst q_{n+2}>\frac{q_{n}}{\g}\su{(\xn)-\frac{q_{n}q_{n+1}}{ q_{n+2}^{\t+1}}-\frac{2q_{n}q_{n+1}}{q_{n+2}^{\t-1}}} \\ 
\end{array}\right.
\ \\\eeq{}
Because of $\a\in\Dgt$, the first of the two conditions is definitively verified, in fact, for $n$ big enough:

$$\frac{q_{n}q_{n+1}}{ q_{n+2}^{\t+1}}+\frac{2q_{n}q_{n+1}}{q_{n+2}^{\t-1}}<
\frac{1}{\a_{n+2}q_{n}^{\t-1}}<\xn$$ 
So, from the fact that $q_{n+2}=a_{n+2}q_{n+1}+q_{n}$ we are done. \qed

\nl
\lem{}
\label{l15}
Let $\t>\frac{3+\sqrt{17}}{2}$. For almost all $\g\in(0,\frac{1}{2})$, if $\a=[a_{0},a_{1},...]\in \Dgt$ ,for $n$ even big enough: (\ref{*}) is true if and only if (\ref{**}) is true.
\elem{}

\proof
If (\ref{**}) is true, then trivially (\ref{*}) is true. 
So we have to show that for almost all $\g\in(0,\su{2})$ and for all $\a\in \Dgt$ (with
$\t>\frac{3+\sqrt{17}}{2}$) holds the converse. So, suppose by contradiction that exists $A\subset \left(C_{1}, C_{2}\right)$, with $0<C_{1}<C_{2}<\su{2}$, $\m(A)>0$ such that, for all $\g\in A$ there exists $\a\in \Dgt$ that satisfies (\ref{*}) but not (\ref{**}) for infinitely many $n$ even. By Lemma \ref{l11} and Lemma \ref{l14} it follows that for all $\g$ in $A$ there exists $\a\in\Dgt$ such that for infinitely many $n$ even:
$$\frac{q_{n}}{\g q_{n+1}}\su{(\xn)-\frac{q_{n}q_{n+1}}{ q_{n+2}^{\t+1}}-\frac{2q_{n}q_{n+1}}{q_{n+2}^{\t-1}}}-\frac{q_{n}}{q_{n+1}}\geq a_{n+2}>\frac{q_{n}}{\g q_{n+1}}\su{(\xn)-\frac{q_{n}q_{n+1}}{ q_{n+2}^{\t+1}}}-\frac{q_{n}}{q_{n+1}},$$
and by Remark \ref{r16} it follows that, for this $n$: $$q_{n+1}\sim \frac{q_{n}^{\t}}{\g}.$$ So, for $n$ big enough such that  (\ref{*}) holds but (\ref{**}) doesn't hold we get:
$$\frac{q_{n}^{\t}}{C_{2}}<q_{n+1}<\frac{q_{n}^{\t}}{C_{1}}.$$
Moreover:
$$a_{n+2}>\frac{q_{n}}{\g q_{n+1}}\su{(\xn)-\frac{q_{n}q_{n+1}}{ q_{n+2}^{\t+1}}}-\frac{q_{n}}{q_{n+1}} \iff$$
$$\frac{a_{n+2}q_{n+1}}{q_{n}}+1=\frac{q_{n+2}}{q_{n}}>\su{1-\frac{\g q_{n+1}}{q_{n}^{\t}}-\frac{\g q_{n}q_{n+1}}{q_{n+2}^{\t+1}}}\iff$$ $$
1-\frac{\g q_{n+1}}{q_{n}^{\t}}-\frac{\g q_{n}q_{n+1}}{q_{n+2}^{\t+1}}>\frac{q_{n}}{q_{n+2}}\iff$$
$$\g<\frac{1-\frac{q_{n}}{q_{n+2}}}{\frac{q_{n+1}}{q_{n}^{\t}}+\frac{q_{n}q_{n+1}}{q_{n+2}^{\t+1}}}$$
In a similar way:
$$\frac{q_{n}}{\g q_{n+1}}\su{(\xn)-\frac{q_{n}q_{n+1}}{ q_{n+2}^{\t+1}}-\frac{2q_{n}q_{n+1}}{ q_{n+2}^{\t-1}}}-\frac{q_{n}}{q_{n+1}}\geq a_{n+2}\iff$$
$$\g\geq\frac{1-\frac{q_{n}}{q_{n+2}}}{\frac{q_{n+1}}{q_{n}^{\t}}+\frac{q_{n}q_{n+1}}{q_{n+2}^{\t+1}}+\frac{2q_{n}q_{n+1}}{q_{n+2}^{\t-1}}}.$$ 
Thus:
$$\frac{1-\frac{q_{n}}{q_{n+2}}}{\frac{q_{n+1}}{q_{n}^{\t}}+\frac{q_{n}q_{n+1}}{q_{n+2}^{\t+1}}+\frac{2q_{n}q_{n+1}}{q_{n+2}^{\t-1}}}\leq \g<\frac{1-\frac{q_{n}}{q_{n+2}}}{\frac{q_{n+1}}{q_{n}^{\t}}+\frac{q_{n}q_{n+1}}{q_{n+2}^{\t+1}}}$$
for infinitely many $n$ even, so for all $\g\in A$ there exist infinitely many $q\in\Bbb{N}$ such that:
$$\frac{1-\frac{q}{Np+q}}{\frac{p}{q^{\t}}+\frac{q p}{(Np+q)^{\t+1}}+\frac{2q p}{ (N p+q)^{\t-1}}}\leq \g<\frac{1-\frac{q}{Np+q}}{\frac{p}{q^{\t}}+\frac{q p}{(Np+q)^{\t+1}}}$$
for some $N\in\Bbb{N}$ and some $\frac{q^{\t}}{C_{2}}<p<\frac{q^{\t}}{C_{1}}$. So for all $M\in \Bbb{N}$:
$$A\subset \bigcup_{q>M}\bigcup_{\frac{q^{\t}}{C_{2}}<p<\frac{q^{\t}}{C_{1}}}\bigcup_{N>0}\left(\frac{1-\frac{q}{Np+q}}{\frac{p}{q^{\t}}+\frac{q p}{(Np+q)^{\t+1}}+\frac{2q p}{(N p+q)^{\t-1}}},\frac{1-\frac{q}{Np+q}}{\frac{p}{q^{\t}}+\frac{q p}{(Np+q)^{\t+1}}}\right),$$
moreover:
$$\frac{1-\frac{q}{Np+q}}{\frac{p}{q^{\t}}+\frac{q p}{(Np+q)^{\t+1}}}-\frac{1-\frac{q}{Np+q}}{\frac{p}{q^{\t}}+\frac{q p}{(Np+q)^{\t+1}}+\frac{2q p}{(N p+q)^{\t-1}}}<$$ $$\frac{2q p}{ (N p+q)^{\t-1}}\left(\su{\frac{p}{q^{\t}}+\frac{q p}{(Np+q)^{\t+1}}}\right)^{2}<\frac{2q C_{2}^{2}}{N^{\t-1}p^{\t-2}}$$
so we obtain:
$$m(A)\leq \sum_{q>M}\sum_{\frac{q^{\t}}{C_{2}}<p<\frac{q^{\t}}{C_{1}}}\sum_{N>0}\frac{2q C_{2}^{2}}{N^{\t-1}p^{\t-2}}<$$
$$\b\sum_{q>M}\frac{q^{\t+1}}{q^{\t^{2}-2\t}}=\b\sum_{q>M}\su{q^{\t^{2}-3\t-1}}$$
for some constant $\b>0$. From the hypothesis ($\t>\frac{3+\sqrt{17}}{2}$) we have that the series converge, so for $M$ that goes to infinity we get that $\m(A)=0$, that contradicts the hypothesis $\m(A)>0$. Thus, for almost all $\g\in(C_{1},C_{2})$ we have that: if (\ref{*}) holds, then (\ref{**}) holds, and from the arbitrariness of $C_{1}, C_{2}$ Lemma \ref{l15} follows.\qed
\\
\nl
\pro{}
\label{p3}
Let $\t>\frac{3+\sqrt{17}}{2}$. For almost every $0<\g<\su{2}$:
if $\a\in\Dgt$, $\np$ are the convergents to $\a$, $\a-\np>\qng$, and definitively:
$$\np+\qng<\nppp-\qnnng,$$
then $\a$ is an accumulation point of $\Dgt$ and in particular, for $n$ even big enough:
$$\m\left(\Dgt\cap \left(\np,\nppp\right)\right)>0$$
\epro{}
\proof
By Lemma \ref{l8} it follows that $\exists N_{1}\in\Bbb{N}$ such that for $n>N_{1}$ even:
$$\pq\not\in I_{n}\implica \pq+\frac{\g}{q^{\t+1}},\pq-\frac{\g}{q^{\t+1}}\not\in A_{n},$$
and by Lemma \ref{l15} for almost all $\g\in(0,\su{2})$:
$$\np+\qng<\nppp-\qnnng\implica \np+\qng<\nppp-\qnnng-\frac{2\g}{q_{n+2}^{\t-1}},$$
therefore, up to a set of measure zero we can suppose that $\g$ satisfies this property. Moreover, by Lemma \ref{l9}, for $n$ even big enough, if $\pq\in I_{n},$ $q<q_{n+2}$ then:
$$\frac{p}{q}+\frac{\g}{q^{\t+1}}<\nppp-\qnnng-\frac{2\g}{q_{n+2}^{\t-1}}.$$
So, if we define: $$c_{n}:=\max_{\pq\in [\np,\nppp),q<q_{n+2}} \pq+\frac{\g}{q^{\t+1}},$$
we obtain:
$$c_{n}<\nppp-\frac{2\g}{q_{n+2}^{\t-1}}-\qnnng.$$
By Lemma \ref{l8}, if $n>N_{1}$ is even and $\pq\not\in I_{n}$, then $$\pq+\frac{\g}{q^{\t+1}},\pq-\frac{\g}{q^{\t+1}}\not\in A_{n},$$ so, if $$\pq<\np \implica \pq+\frac{\g}{q^{\t+1}}<\np+\qng\leq c_{n},$$ while for $\pq>\nppp$ we get $q>q_{n+2}$, so: $$\pq-\frac{\g}{q^{\t+1}}>\nppp-\qnnng,$$
but from: $$\b\in \Dgt^{c}\iff \exists \pq\in(0,1):\b\in \left(\pq-\frac{\g}{q^{\t+1}},\pq+\frac{\g}{q^{\t+1}}\right)$$ we get that for $n>N_{1}$ even, holds:
$$\m\left(\Dgt^{c}\cap I_{n}\right)\leq \m\left(\bigcup_{\pq\in [\np,\nppp),q<q_{n+2}} \left(\pq-\frac{\g}{q^{\t+1}},\pq+\frac{\g}{q^{\t+1}}\right)\cap I_{n} \right)$$ $$+\m\left(\bigcup_{\pq \in I_{n},q\geq q_{n+2}} \left(\pq-\frac{\g}{q^{\t+1}}, \pq+\frac{\g}{q^{\t+1}}\right)\right)+\m\left(\nppp-\qnnng,\nppp\right).$$
So by Lemma \ref{l10}:
$$\m(\Dgt^{c}\cap I_{n})\leq c_{n}-\np+\frac{2\g}{q_{n+2}^{\t-1}}+\qnnng<\m(I_{n})=\nppp-\np \iff$$
$$c_{n}<\nppp-\qnnng-\frac{2\g}{q_{n+2}^{\t-1}},$$
that follows from the definition of $c_{n}$. \qed
So, given $\t>3$, for almost all $\g>0$:
if $\a\in\Dgt$ is not an isolated point of the first type and definitively the intervals centered in the convergents have an empty intersection, then $\a$ is an accumulation point in $\Dgt$.
The second step is to show that: if $\t>3$, $\g>0$, $\a\in\Dgt$ but $\a$ is not an isolated point of the first type and $\t>\t(\a)$, then $\a$ is an accumulation point in $\Dgt$.

\nl
\lem{}
\label{l16}
Let $\t>3$. For almost all $\g\in(0,\su{2})$: given $\a\in \Dgt$, if for infinitely many $n$ even:
$$\np+\qng>\nppp-\qnnng,$$
then there exists $C>0$ such that for this $n$:
$$a_{n+2}\leq C q_{n}^{2+\e},$$
with $\e>0$ arbitrarily small.
\elem{}
\proof
By Lemma \ref{l11} it follows that, given $\a\in \Dgt$ that satisfies the hypothesis of Lemma \ref{l16}, for $n$ even big enough:
 $$a_{n+2}\leq\frac{q_{n}}{\g q_{n+1}}\su{(\xn)-\frac{q_{n}q_{n+1}}{ q_{n+2}^{\t+1}}}-\frac{q_{n}}{q_{n+1}},$$
 so, up to a set of measure zero, by Lemma \ref{l12} we can suppose that there exist $\e>0, C>0$ such that $\su{\g}\in B_{C,\t+1+\e}$ with $\t+1+\e<\t^{2}-1$, from which it follows that:
$$ \frac{q_{n}}{\g q_{n+1}}\su{(\xn)-\frac{q_{n}q_{n+1}}{ q_{n+2}^{\t+1}}}-\frac{q_{n}}{q_{n+1}}\leq \frac{q_{n}}{\g q_{n+1}}\su{\frac{C}{q_{n}^{\t+1+\e}}-\frac{q_{n}q_{n+1}}{ q_{n+2}^{\t+1}}}-\frac{q_{n}}{q_{n+1}},$$
moreover, by Remark \ref{r16} it follows that $q_{n+1}\sim \frac{q_{n}^{\t}}{\g}$,
from which we obtain:
$$\frac{q_{n}q_{n+1}}{ q_{n+2}^{\t+1}}<\frac{q_{n}}{q_{n+1}^{\t}}\sim \frac{\g^{\t}}{q_{n}^{\t^{2}-1}},$$
so, if $n$ is big enough, by $\t+1+\e<\t^{2}-1$ we have:
$$\frac{C}{q_{n}^{\t+1+\e}}-\frac{q_{n}q_{n+1}}{ q_{n+2}^{\t+1}}>\frac{C}{2q_{n}^{\t+1+\e}}.$$
So we obtain:
$$a_{n+2}<\frac{q_{n}}{q_{n+1}}\frac{2 q_{n}^{\t+1+\e}}{C}\sim \frac{2\g}{ C}q_{n}^{2+\e}<\frac{4\g}{ C}q_{n}^{2+\e}=C' q_{n}^{2+\e}$$
definitively, from which we get Lemma \ref{l16}. \qed
\\
\nl
\lem{}
\label{l17}
Let $\t>\frac{3+\sqrt{17}}{2}, \g>0$, $\a\in\Dgt$. If for infinitely many $m$ even, for $n<m$ even holds:
\beq{a}
\np+\qng<\qm-\qmg-\frac{2\g}{q_{m}^{\t-1}},
\eeq
and $\a-\np>\qng$ for all $n$ even, then $\a$ is in ${\mathcal{A}}(\Dgt)$.
\elem{}
\proof
Let $\np<\pq<\nppp$ with $n$ even and $n<m-2$, for $\frac{q_{n+2}}{2}\leq q$:
$$\pq+\frac{\g}{q^{\t+1}}<\nppp+\qnnng$$
is definitively true, while for $q<\frac{q_{n+2}}{2}$:
$$\nppp-\pq=\nppp-\a+\a-\pq>\frac{1}{2q^{2}}-\frac{1}{q_{n+2}^{2}}>\frac{\g}{q^{\t+1}}-\qnnng\iff$$
$$\frac{1}{2q^{2}}-\frac{\g}{q^{\t+1}}>\frac{1}{q_{n+2}^{2}}-\qnnng,$$
that is true for $q$ big enough, so $\exists T\in\Bbb{N}$ such that the inequality is verified for $q\geq T$ (from the fact that $G(x):=\frac{1}{2x^{2}}-\frac{\g}{x^{\t+1}}$ is definitively decreasing and $\t>3>1$).
From the hypothesis that $\a-\np>\qng$ for all $n$ even:  $$v:=\max_{\frac{p}{q}<\a,q\leq T}\pq+\frac{\g}{q^{\t+1}}<\a,$$
so there exists $T_{1}\in\Bbb{N}$ such that for $n>T_{1}$:
$$\np+\qng>v.$$
By Lemma \ref{l9}, for $m$ big enough, $\frac{p}{q}\in I_{n}$, with $n<m-2$ even:
$$\pq+\frac{\g}{q^{\t+1}}\leq\max\left\{\nppp+\qnnng,v\right\}\leq \frac{p_{m-2}}{q_{m-2}}+\frac{\g}{q_{m-2}^{\t+1}}<\qm-\qmg-\frac{2\g}{q_{m}^{\t+1}},$$
while by Lemma \ref{l9}, for $m$ big enough:
$$\pq\in I_{m-2},q<q_{m-2}\implica \pq+\frac{\g}{q^{\t+1}
}<\qm-\qmg-\frac{2\g}{q_{m}^{\t-1}},$$
so if we define:
$$c_{m}:=\max\left\{\max_{\frac{p}{q}\in I_{m-2}, q<q_{m}} \left(\pq+\frac{\g}{q^{\t+1}}\right), \max_{\frac{p}{q}\leq \frac{p_{m-2}}{q_{m-2}}} \left(\pq+\frac{\g}{q^{\t+1}}\right)\right\},$$
for $m$ even big enough:
$$c_{m}<\qm-\qmg-\frac{2\g}{q_{m}^{\t-1}}.$$
Moreover, by Lemma \ref{l10}, from $\t>3>2$, for $m$ even big enough:
$$\m\left(\bigcup_{\frac{p}{q}\in I_{m-2}, q\geq q_{m}}\left(\frac{p}{q}-\frac{\g}{q^{\t+1}}, \frac{p}{q}+\frac{\g}{q^{\t+1}}\right)\right)<\frac{2\g}{q_{m}^{\t-1}}.$$
Finally, if $\pq>\frac{p_{m}}{q_{m}}$, by the properties of continued fractions we obtain $q>q_{m}$, so $\frac{p}{q}-\frac{\g}{q^{\t+1}}>\qm-\qmg$.
Thus:
$$\m\left(\Dgt^{c} \cap\left(\frac{p_{m-2}}{q_{m-2}},\qm-\qmg\right)\right)<c_{m}-\frac{p_{m-2}}{q_{m-2}}+\frac{2\g}{q_{m}^{\t-1}}$$ $$<\qm-\frac{p_{m-2}}{q_{m-2}}-\qmg=\m(\frac{p_{m-2}}{q_{m-2}},\qm-\qmg),$$
then $$\Dgt \cap\left(\frac{p_{m-2}}{q_{m-2}},\qm-\qmg\right)\not =\emptyset,$$ and from the fact that this holds for infinitely many $m$ even, then $\a$ is an accumulation point of $\Dgt$. \qed
\\
\nl
\rem{}
\label{r17}
Let $\t>\frac{\sqrt{17}+3}{2},\g>0$, $\a\in\Dgt$, if $\a\in \iigt$ or $\iiigt$, then $\t(\a)=\t$. In fact if this doesn't hold, from $\a\not\in \igt$ we get that for all $n$ even or for all $n$ odd:
$$\left|\a-\np\right|>\frac{\g}{q_{n}^{\t+1}}.$$
Suppose for example that this property holds for all $n$ even. If on the contrary $\t(\a)<\t$, by  Remark \ref{r16}, the hypothesis of Proposition \ref{p3} are satisfied, so $\a\in {\mathcal{A}}(\Dgt)$, contradiction.   
\erem{}

\nl
\cor{}
\label{c4}
If $\t>\frac{3+\sqrt{17}}{2}$: $$\m\left(\left\{\g>0: \iigt\not=\emptyset\right\}\right)=0.$$
\ecor{}
\proof
Observe that, if $\a\in \iigt$, then there exists $n\in\Bbb{N}$ such that $$\left|\a-\np\right|=\qng.$$ Suppose for example that $n$ is even, thus:
$$\a=\np+\qng.$$
Moreover, for almost all $\g\in(0,\su{2})$:
$$\t\left(\pq+\frac{\g}{q^{\t+1}}\right)=\t\left(\frac{\g}{q^{\t+1}}\right)=1.$$
Taking the union on all the $\pq$ we obtain that for almost all $\g\in(0,\su{2})$ and for all $\pq\in\Bbb{Q}$, $$\t\left(\pq+\frac{\g}{q^{\t+1}}\right)=1.$$
So Corollary \ref{c4} follows by Remark \ref{r17}.\qed

\nl
It remains the last one step, in which we get the Theorem.
\lem{}
\label{l18}
Let $\t>3$. For almost all $\g>0$, if $\a\in {\II}(\Dgt)$, there exists $N\in \Bbb{N}$ such that, for all $m>N$ even there is some $n<m$ even with:
$$\np+\qng \geq \qm-\qmg-\frac{2\g}{q_{m}^{\t-1}}$$.
\elem{}
\proof
By Corollary \ref{c4} and Remark (j) it follow that, up to a set of measure zero, we can suppose that $\igt=\iigt=\emptyset$, so observe that if the Lemma were not true, it would exist $\a\in \iigt$ with the even convergents that satisfy the hypothesis of Lemma \ref{l17}, that implies $\a\in {\mathcal{A}}(\Dgt)$, contradiction. \qed
\\
\nl
{\bf Theorem}
Let $\t>\frac{3+\sqrt{17}}{2}$. Then, for almost all $\g>0$ $\Dgt$ is a Cantor set.

\proof
By Corollary \ref{c4} and Remark (j) it follows that, up to a set of measure zero, we can suppose that $\igt=\iigt=\emptyset$. Suppose by contradiction that the statement doesn't hold, and take $0<C_{1}<C_{2}$ such that:
$$\m\left(\left\{C_{1}<\g<C_{2}: {\II}(\Dgt)\not=\emptyset\right\}\right)>0,$$
and define $A:=\{C_{1}<\g<C_{2}: {\II} (\Dgt)\not=\emptyset\}$. By Lemma \ref{l18}, for almost all $\g>0$ there exists $\a\in {\II}(\Dgt)$ and there exists $N\in\Bbb{N}$ such that for all $m>N$ even, there is some $n<m$ even, with:
$$\np+\qng\geq \qm-\qmg-\frac{2\g}{q_{m}^{\t-1}}.$$
Now we want to show that, for almost all chosen of $\g\in A$ we have: 
$$\limsup \frac{q_{2k+2}}{q_{2k+1}^{\t}}<\frac{1}{\g}.$$
In fact if it doesn't hold, by Remark \ref{r16} we get that for infinitely many $m$ even:
$$ q_{m}\sim \frac{q_{m-1}^{\t}}{\g},$$
and for $m>N$ exists $n<m$ even, with:
$$\np+\qng\geq \qm-\qmg-\frac{2\g}{q_{m}^{\t-1}}$$
By Lemma \ref{l15}, up to a set of measure zero in $A$:
$$\np+\qng\geq \qm-\qmg-\frac{2\g}{q_{m}^{\t-1}}\iff \np+\qng\geq \qm-\qmg.$$
By the properties of convergents:
$$\a-\qm<\frac{1}{q_{m}^{2}},$$
from which we get:
$$\frac{1}{q_{m}^{2}}>\a-\np-\qng-\qmg.$$
Moreover:
$$\a-\np=\su{q_{n}(q_{n+1}+\frac{\a_{n+2}}{q_{n}})},$$
so:
$$\frac{1}{q_{m}^{2}}>\su{q_{n}(q_{n+1}+\frac{\a_{n+2}}{q_{n}})}-\qng-\qmg$$
For $m$ big enough:
$$\frac{1}{q_{m}^{2}}+\qmg<\frac{2}{q_{m}^{2}},$$
so:
$$\frac{2}{q_{m}^{2}}>\su{q_{n}(q_{n+1}+\frac{\a_{n+2}}{q_{n}})}-\qng\iff$$
$$\g>\frac{q_{n}^{\t}}{q_{n+1}+\frac{\a_{n+2}}{q_{n}}}-\frac{2q_{n}^{\t+1}}{q_{m}^{2}},$$ 
moreover:
$$\g\leq\frac{q_{n}^{\t}}{q_{n+1}+\frac{\a_{n+2}}{q_{n}}}.$$
So we obtain:
$$\frac{q_{n}^{\t}}{q_{n+1}+\frac{\a_{n+2}}{q_{n}}}-\frac{2q_{n}^{\t+1}}{q_{m}^{2}}<\g\leq\frac{q_{n}^{\t}}{q_{n+1}+\frac{\a_{n+2}}{q_{n}}}$$
From
$$\np+\qng\geq \qm-\qmg-\frac{2\g}{q_{m}^{\t-1}},$$
we get:
$$\np+\qng\geq \nppp-\qnnng,$$
moreover, from $\a-\np>\qng$ for all $n$ even, when $m$ increase, also $n$ increase, and by the last inequality and Remark \ref{r16} we get that $q_{n+1}\sim \frac{q_{n}^{\t}}{\g}$. So
$$q_{m}\sim \frac{q_{m-1}^{\t}}{\g}\geq \frac{q_{n+1}^{\t}}{\g}\sim \frac{q_{n}^{\t^{2}}}{\g^{\t}}\geq\frac{q_{n}^{\t^{2}}}{C_{2}^{\t}}.$$
So we obtain:
$$\frac{q_{n}^{\t}}{q_{n+1}+\frac{\a_{n+2}}{q_{n}}}-\frac{C}{q_{n}^{2\t^{2}-\t-1}}<\g\leq\frac{q_{n}^{\t}}{q_{n+1}+\frac{\a_{n+2}}{q_{n}}}$$
with a constant $C>0$. By Lemma \ref{l16}, up to a set of measure zero, we can suppose that there exists $\e>0$ arbitrarily small such that, for $n$  big enough:
$$a_{n+2}<q_{n}^{2+\e}.$$
So, up to a set of measure zero, we can suppose that for all $\g\in A$, there exists infinitely many $q>0$, $\frac{q^{\t}}{2 C_{2}}<p<\frac{2}{C_{1}q^{\t}}$, $N<q^{2+\e}$ such that:
$$\frac{q^{\t}}{p+\frac{N}{q}}-\frac{C}{q^{2\t^{2}-\t-1}}<\g\leq\frac{q^{\t}}{q+\frac{N}{q}}.$$ So, for all $M\in\Bbb{N}$:
$$A\subset \bigcup_{q>M}\bigcup_{\frac{q^{\t}}{2C_{2}}<p<\frac{2q^{\t}}{C_{1}}}\bigcup_{N<q^{2+\e}} \left(\frac{q^{\t}}{p+\frac{N}{q}}-\frac{C}{q^{2\t^{2}-\t-1}}, \frac{q^{\t}}{q+\frac{N}{q}}\right),$$
Thus:
$$\m(A)<\sum_{q>M}\sum_{\frac{q^{\t}}{2C_{2}}<p<\frac{2q^{\t}}{C_{1}}}\sum_{N<q^{2+\e}} \frac{C}{q^{2\t^{2}-\t-1}}$$
$$<\b \sum_{q>M}\su{q^{2\t^{2}-2\t-3-\e}}$$
with some constant $\b>0$. Because of $\t>\frac{3+\sqrt{17}}{2}$, for $\e$ small enough the series converge, so for $M$ that tends to infinity we obtain $\m(A)=0$, contradiction.
So we have proved that:
$$\limsup \frac{q_{2k+2}}{q_{2k+1}^{\t}}<\frac{1}{\g}.$$
But, by Remark \ref{r16} and Proposition \ref{p3} (used with $n$ odd) we have that $\a\in {\mathcal{A}}(\Dgt)$, contradiction. So $\m(A)=0$.\qed
\nl
The estimate $\t>\frac{3+\sqrt{17}}{2}$ can be improved putting a better inequality in Lemma 5. Probably  the Proposition holds also with $\t>3$.
\section{Questions}
\begin{itemize}
    \item By \cite{17} we konw that, for some choice of $\g,\t$, $\igt\not=\emptyset$. What about $\iigt,\iiigt$?
    \item What is the best $\t>1$ such that the result holds?
    \item Is it true that, for all $\t\geq1$ there exists $\g_\t\in(0,\su{2})$ such that $\Dgt$ is a Cantor set for almost all $\g\in(0,\g_\t)$?
\end{itemize}
\subsubsection*{Acknowledgement}
I am very grateful to Prof. Luigi Chierchia for his suggestions, remarks, for his special support and for encouraging me to complete this work.

\end{document}